\documentclass[10pt,a4paper]{gOMS2e}
\pdfoutput=1 

\usepackage{geometry}

\usepackage{bm}
\usepackage{color}
\usepackage{courier}         %
\usepackage{makeidx}         %
\usepackage{graphicx}        %
\usepackage{listings}        %
\usepackage{multicol}        %
\usepackage[bottom]{footmisc}%

\usepackage{amsmath}         %
\usepackage{amssymb}         %
\usepackage{blkarray}        %
\usepackage{enumerate}       %
\usepackage{ifthen}          %
\usepackage{xspace}          %
\usepackage{mathtools}       %
\usepackage{url}             %
\usepackage{wrapfig}         %
\newcommand\nc{\newcommand}
\nc\rnc{\renewcommand}

\newcommand\bc{\begin{center}}
\newcommand\ec{\end{center}}
\newcommand\bn{\begin{enumerate}}
\newcommand\en{\end{enumerate}}
\newcommand\TW{\textwidth}

\makeatletter
\def\rf#1{(\@rf#1,.)}
\def\@rf#1,{\ref{eq:#1}\@ifnextchar . {\@endrf}{, \@rf}}
\def\@endrf.{}
\makeatother

\newcommand\exref[1]{Example~\ref{ex:#1}\xspace}
\newcommand\fgref[1]{Figure~\ref{fg:#1}\xspace}

\newcommand\scref[1]{Section~\ref{sc:#1}\xspace}

\newcommand\tbref[1]{Table~\ref{tb:#1}\xspace}

\newcommand\ssrf[1]{\S\ref{ss:#1}\xspace}%

\newcommand\mc[3]{\multicolumn{#1}{#2}{#3}}
\renewcommand\_[1]{\mathbf{#1}}

\newcommand\ninf{{-\!\infty}}

\renewcommand{\d}{\mathrm{d}}

\newcommand\lam{\lambda}
\newcommand\M{\mathcal{M}}

\newcommand\sij{\sigma_{ij}}

\newcommand\set[2]{\{\,#1\mid\mbox{#2}\,\}}
\newcommand\dbd[2]{\frac{\partial #1}{\partial #2}} %
\newcommand\dbdt[2]{\partial #1/\partial #2} %

\newcommand\matlab{{\sc Matlab}\xspace}

\newcommand\SA{structural analysis\xspace}
\newcommand\Smethod{$\Sigma$-method\xspace}
\newcommand\sigmx{signature matrix\xspace}

\newcommand\sysJ{system Jacobian\xspace}

\nc\HI[1]{{\color{red}#1}}

\nc\Apd{\color[rgb]{0,.6,0}}
\nc\Av{\color[rgb]{.2,.5,.7}}
\nc\Gd{\color[rgb]{.4,.2,1}}

\newtheorem{lemma}{Lemma}[section]
\newtheorem{examp}{{Example}}[section]
\newenvironment{example}[1]{\begin{examp}[#1]\rm}{\qquad\qed\end{examp}}

\def\@#1@{\texttt{\small #1}}
\nc\kc{{k_\textup{c}}}
\nc\kd{{k_\textup{d}}}
\nc\Kmax{\kappa_\textup{max}}
\nc\bmx[1]{\begin{bmatrix}#1\end{bmatrix}}
\nc\sbmx[1]{\left[\begin{smallmatrix}#1\end{smallmatrix}\right]}
\nc\pmx[1]{\begin{pmatrix}#1\end{pmatrix}}
\rnc\.[2]{\ifthenelse{#1=1}{\dot{#2}}{\ifthenelse{#1=2}{\ddot{#2}}{\ifthenelse{#1=3}{\dddot{#2}}{??}}}}
\rnc\_[1]{\mathbf{#1}}
\rnc\d{\textup{d}}
\nc\bx{\_x}
\rnc\x{\times}
\nc\bdelta{\bm{\delta}} %
\nc\der[3]{#1_{#2}^{(#3)}} %
\nc\oo{\infty}

\nc\AD{algorithmic differentiation\xspace}
\nc\cpp{{\tt C++}\xspace}
\nc\CM{centre of mass\xspace}
\nc\daets{{\sc Daets}\xspace}
\nc\DOF{\textsc{dof}\xspace}
\nc\MS{Mattsson \& S\"oderlind\xspace}
\nc\fadbadpp{{\tt FADBAD++}\xspace}
\nc\found{found\xspace}
\nc\DAE{differential-algebraic equation\xspace}
\nc\DAEs{\DAE{s}\xspace}
\nc\EBM{equation-based modeling\xspace}
\nc\genpend{\textsc{GenericPend}\xspace}
\nc\gproms{gPROMS\xspace}
\nc\wrt{with respect to\xspace}
\nc\sa{structural analysis\xspace}

\nc\li[1]{\lstinline[basicstyle=\tt]{#1}}
\title{How AD Can Help Solve Differential-Algebraic Equations}
\thanks{\copyright 2017 Pryce, Nedialkov, Tan, Li, \qquad March, 21 2017}

\author{
\name{John D Pryce$^a{}^{\ast}$
and Nedialko S Nedialkov$^b$
and Guangning Tan$^c$
and Xiao Li$^d$
}
\affil{$^a$Cardiff University School of Mathematics, Senghennydd Rd, Cardiff, CF24 4AG, Wales;\\
$^b$McMaster University Department of Computing \& Software, Hamilton, L8S 4K1, Canada.\\
$^c$Massachusetts Institute of Technology, Process Systems Engineering Laboratory,  Boston, USA.\\
$^d$McMaster University School of Computational Science and Engineering, Hamilton, L8S 4K1, Canada.
}
}%

\date{April 2016}

\usepackage{fancyheadings}
\begin{document}

\maketitle
\begin{abstract}
  A characteristic feature of differential-algebraic equations is that one needs to find derivatives of some of their equations with respect to time, as part of so called index reduction or regularisation, to prepare them for numerical solution.
  This is often done with the help of a computer algebra system.
  We show in two significant cases that it can be done efficiently by pure algorithmic differentiation.
  The first is the Dummy Derivatives method; here we give a mainly theoretical description, with tutorial examples.
  The second is the solution of a mechanical system directly from its Lagrangian formulation. Here we outline the theory and show several non-trivial examples of using the ``Lagrangian facility'' of the Nedialkov-Pryce initial-value solver DAETS, namely: a spring-mass-multipendulum system; a prescribed-trajectory control problem; and long-time integration of a model of the outer planets of the solar system, taken from the DETEST testing package for ODE solvers.
\end{abstract}

\section{Introduction}\label{sc:intro}

\subsection{DAE formulation and basic ideas}\label{ss:basicideas}
In industrial engineering, the modeling of systems to simulate their time evolution is increasingly done by methods that lead to a {\DAE} (DAE) system as the underlying mathematical form.
Such DAEs often come from equation-based modeling (EBM), which describes system components by the basic physical laws they obey and supports ``multiphysics'' models that combine several scientific disciplines, as for instance mechanical, electrical, chemical, and thermodynamic behaviour in a car engine.

Facilities created to support EBM include \gproms, which is both a language and a graphical modeling environment (GME) built on it; the Modelica language and GMEs such as OpenModelica, Dymola and MapleSim that are built on it. Simulink, built on \matlab, is a GME of similar scope but less in tune with the general DAE concept.

A DAE is just a set of $n$ equations connecting a vector $\_x=\_x(t)$ of $n$ state variables $x_1,\ldots,x_n$ and some derivatives of them \wrt time $t$.
One can always reduce it to a first order form $\_F(t,\_x,\.1{\_x})=0$---as accepted by the DASSL solver and its relatives \cite{BrenanCampbelPetzold,sundials}---in the same way as one does for an ODE system.
Here $\.1{\_x}$ means $\d\_x/\d t$.
However we use a more flexible form allowing arbitrary higher derivatives:
\begin{align}\label{eq:maineq}
  f_i(\, t,\, \textrm{the $x_j$ and derivatives of them}\,) = 0, 
  \quad i=1,\ldots,n.
\end{align}
This often lets one formulate problems to our \daets initial-value code \cite{nedialkov2008solving,nedialkov2008daets} more concisely, e.g.\ Lagrange's equations for a mechanical system with $n_q$ coordinates and $n_c$ constraints need $n_q+n_c$ variables, compared to $2n_q+n_c$ in first-order form.

\subsection{Aim}

In general, differentiating some of the DAE's equations $f_i=0$ \wrt $t$ is an essential step in solving a DAE.
This article is about two significant and rather different uses of this. 
The first is the widely used {\em dummy derivatives} method of \MS \cite{Matt93a} that prepares a higher-index DAE for numerical solution by a classical index-1 DAE code, or by an explicit ODE code such as a Runge--Kutta method.

The second is the task of solving a, possibly constrained, mechanical system directly from a Lagrangian formulation.
Conceptually it has a ``two-phase'' aspect. The motion of the system is defined by a {\em Lagrangian function} $L(q,\.1q)$ where $q=(q_1,\ldots,q_n)$ is a vector of generalised coordinates, together with $m$ constraints $C(q) = \bigl(C_1(q),\ldots,C_m(q)\bigr)=0$.
To set up (phase 1) the equations of motion from $L$ and $C$ involves partial differentiation $\dbdt{}{q}$ and $\dbdt{}{\.1q}$ applied to $L$ and $C$, as well as straight differentiation $\d/\d t$.
When $m>0$ the result is a DAE, which must (phase 2) be readied for numerical solution.

\smallskip
Each of these use cases at first sight seems to need symbolic differentiation, e.g.\ in a computer algebra system.
We show pure AD suffices in either case. In the second case, which we devote more space to, AD can even do the two ``phases'' simultaneously, and gives a simple and elegant user interface and an efficient numerical solution process.

\subsection{Structural analysis}\label{ss:SA}

In an ODE $\.1{\_x}=\_f(t,\_x)$, causality is obvious: in differential language, it explicitly specifies the state $\_x+\d\_x$ at the next instant $t+\d t$ to be $\_x + \_f(t,\_x)\d t$.

In a DAE, causality is not obvious.
For instance, these size 2 DAEs are quite different, where $u(t)$ is a given driving function:
\begin{align}
                   x_1-u(t) &= 0,\quad x_1-\.1x_2 = 0, \label{eq:ex1a} \\
  \text{and}\quad  x_2-u(t) &= 0,\quad x_1-\.1x_2 = 0 \label{eq:ex1b}
\end{align}

To solve \rf{ex1a}, make $\.1x_2$ the subject of its second equation ($x_1$ causes $x_2$) and integrate the result; it is really an ODE, with one degree of freedom.
To solve \rf{ex1b}, make $x_1$ the subject of its second equation ($x_2$ causes $x_1$) and differentiate. DAE \rf{ex1b} has no degrees of freedom---it has the unique solution $x_1=\.1u(t)$, $x_2=u(t)$ and does not look like an ODE at all; such behaviour is common in control problems.

\smallskip
A solvable DAE has a chain of causality that must be found in order to prepare for numerical solution.
Knowing which equations $f_i=0$ to differentiate, and how often, is crucial to finding this causal chain.
When correctly done, the original DAE augmented by the differentiated equations can be solved to produce an ODE in some (possibly not all) of the original variables---the {\em ODE part}.
Once this ODE is solved, the remaining variables forming the {\em algebraic part} can be found by algebraic manipulations combined with differentiations.

Let $c_i$ be the number of differentiations of equation $i$ needed by the ``most economical'' way of doing this.
For reasons to do with the Taylor series method used by \daets we call them the {\em equation-offsets}.

For instance the equations of \rf{ex1a} do not need differentiating: $(c_1,c_2)=(0,0)$. We solve to produce the ODE part $\.1x_2=u(t)$ in just $x_2$. 
By contrast, \rf{ex1b} has $(c_1,c_2)=(1,0)$ meaning the first equation must be differentiated, after which we solve to get $x_1=\.1u(t)$, $x_2=u(t)$. The ODE part is empty.

In the DAE \rf{ex1a}, it happens we can solve for the algebraic variable $x_1$ to get $x_1=u(t)$, indep\-endently of solving the ODE, but this need not be so: if we change it to
\begin{align}
     x_1-x_2-u(t) &= 0,\quad x_1-\.1x_2 = 0, \label{eq:ex1c}
\end{align}
then the ODE part, namely $\.1x_2-x_2-u(t)=0$, must be solved before we know $x_1$.

\smallskip
Unlike a well-behaved ODE $\.1{\_x}=\_f(t,\_x)$, which has a solution path through each point of the region $R$ of $(t,\_x)$ space where it is defined, the union of a typical DAE's solution paths is a proper subset of $R$, the {\em consistent manifold} $\M$ or set of {\em consistent points}.
The dimension of its intersection with any time $t=t_0$ is \DOF, the number of {\em degrees of freedom}, equivalently the size of its ODE part (here assumed independent of $t_0$).

\smallskip
The {\em index} of a DAE used in this paper is simply
\begin{align}\label{eq:indexdef}
  \nu = \max_i c_i.
\end{align}
The classic {\em differentiation index} $\nu_d$ assigns index $1$ and $2$ to these DAEs respectively.
In summary for the examples above
\[
\begin{tabular}{c|cccccc}
  DAE & ODE part & \DOF & algebraic part & offsets & $\nu$ & $\nu_d$ \\\hline
  eqns \rf{ex1a} & $x_2$ & $1$ & $x_1$ (found independently of ODE part) & $(0,0)$ & $0$ & $1$ \\
  eqns \rf{ex1b} & empty & $0$ & $x_1,x_2$ & $(1,0)$ & $1$ & $2$ \\
  eqns \rf{ex1c} & $x_2$ & $1$ & $x_1$ (found using $x_2$ in ODE part) & $(0,0)$ & $0$ & $1$
\end{tabular}
\]

The {\em\sa} (SA) approach aims to derive a DAE's causal chain by studying its sparsity, namely what derivatives of variables occur in what equations.
The method is: seek a number $c_i$ of times to differentiate the $i$\,th equation that gives a {\em structurally nonsingular} (SNS) set of equations for the resulting highest, $d_j$\,th, derivatives of the $x_j$---then $\_c = (c_1,\ldots,c_n)$, $\_d = (d_1,\ldots,d_n)$ are the vectors of equation-offsets and corresponding variable-offsets.
SNS means one can make a matching of variables to equations---equivalently a {\em transversal}, a set $T$ of $n$ positions $(i,j)$ in an $n\x n$ matrix with just one in each row and in each column---such that derivative $\der xj{d_j}$ occurs in the differentiated equation $\der fi{c_i}=0$ for each $(i,j)\in T$.
There exist unique elementwise smallest non-negative $\_c,\_d$, the {\em canonical offsets}, which we assume chosen henceforth.
They define the ``most economical'' differentiations mentioned above.

An {\em SA-friendly} DAE by definition is one for which these equations are actually (not just structurally) nonsingular at some consistent point, that is, the $n\x n$ {\em system Jacobian}
\begin{align}\label{eq:sysJ}
 \_J = \left(\dbdt{\der fi{c_i}}{\der xj{d_j}}\right)_{i,j=1,\ldots,n}.
\end{align}
is nonsingular there. 
Assuming suitable smoothness of the $f_i$, a unique solution then exists locally through this point, and through any nearby consistent points.

Experience shows most DAEs in practice are SA-friendly. This fact underlies the wide use of the dummy derivatives method, which uses the results of SA and succeeds if and only if the DAE is SA-friendly. The SA can be done by the graph-based Pantelides method \cite{Pant88b}, or the Pryce \Smethod \cite{Pryce2001a} based on the {\em \sigmx} $\Sigma = (\sij)$, where 
\begin{align}\label{eq:sigmx}
  \sij = 
\begin{cases}
  \text{order of highest derivative of $x_j$ in $f_i$} &\text{if $x_j$ occurs in $f_i$}, \\
  -\oo &\text{if not}.
\end{cases}
\end{align}
The methods are equivalent except that the latter handles higher-order DAEs without reduction to first order, while the former as described in \cite{Pant88b} does not.

The DAE (of classical index $3$) derived from a constrained Lagrangian of a mechanical system as in \scref{lagrangian}, is always SA-friendly when posed as an initial value problem.
Posed otherwise, e.g.\ as a prescribed-trajectory control problem, it need not be.
The occurrence of non-SA-friendly but solvable DAEs in applications is studied in \cite{schwarz1998structural,riaza2007qualitative}. For systematic ways of converting such a DAE to an equivalent SA-friendly one see \cite{Tan2015AMMCS}.

\section{Dummy derivatives}\label{sc:DDs}

\subsection{The DDs construction}\label{ss:DDs}

Many numerical methods for higher-index DAEs start with {\em index reduction}: augmenting the DAE by time-derivatives of some of its equations to produce a DAE of larger size and smaller index.
Various index reduction methods have been used that convert the DAE to an ODE with more degrees of freedom than the DAE. Then the solution paths of the DAE are a proper subset of those of the ODE.
This tends to be bad numerically, as errors cause drift from the consistent manifold that is often exponential once it starts.

Dummy derivatives (DDs) by contrast are a systematic way to form an equivalent ODE with {\em exactly as many} DOF as the DAE.
If one regards the DAE as a flow on the consistent manifold $\M$, DDs describes the flow in a local coordinate system for $\M$.
Thus numerical drift can only be within $\M$, where it is less harmful.
However if the path leaves the patch of $\M$ where the coordinate system is nonsingular, one must choose new coordinates. This need for {\em DD switching}, or {\em pivoting}, complicates a numerical algorithm.

\smallskip
The following description of the DDs process is equivalent to that in \cite{Matt93a}.
The set of possible matrix sequences $(\_G_k)$ whenever one selects a state vector, below, is the same in either method, but we find $\_G_k$ from smallest up (each is a sub-matrix of the next), while \cite{Matt93a} finds them in the opposite order.

Assume $c_i$ and $d_j$ are the canonical offsets.
First form the derivatives of each $f_i=0$ up to the $c_i\,$th, forming the augmented system of $N_f=n+\sum_i c_i$ equations:
\begin{align}\label{eq:fider}
  \der fil = 0, \quad l=0,\ldots,c_i,\; i=1,\ldots,n.
\end{align}
Its unknowns are the $N_x=n+\sum_j d_j$ derivatives of the state variables $x_j$ up to the $d_j\,$th.
View them for now as unrelated algebraic unknowns that we call {\em items}, and to emphasise this denote them $x_{jl}$:
\begin{align}\label{eq:xjder}
  x_{jl} \text{ renames } \der xjl, \quad l=0,\ldots,d_j,\; j=1,\ldots,n.
\end{align}
The system has fewer equations than variables by the amount $\sum_j d_j - \sum_i c_i$, which equals the number \DOF of degrees of freedom.
To balance this, DDs finds a number \DOF of items $x_{jl}$ to be {\em state items}, for $(j,l)$ in a suitable set $S$ of index pairs, chosen such that all the other items can locally be solved for as functions of these.
The {\em state vector} $\_x_S$ formed by the state items is the associated local coordinate system of the manifold $\M$.

One requires $l<d_j$ for each $(j,l)\in S$, so that $x_{j,l+1}$ is also an item.
Then the differential relations between each state item and its next higher derivative:
\begin{align}\label{eq:livederiv}
  \.1x_{jl} = x_{j,l+1}
\end{align}
can be interpreted as an ODE system for the state items.

State vector selection---initially or at a DD-switching point---may be done as follows.
The $n\x n$ \sysJ $\_J$ in \rf{sysJ} is nonsingular there.
For $k=\kd,\kd+1,\ldots,-1$ where $\kd$ is minus the largest $d_j$, the ``standard solution scheme'' of the \Smethod constructs sub-matrices $\_J_k$ of $\_J$ by selecting those rows $i$ where $k+c_i\ge0$ and columns $j$ where $k+d_j\ge0$. 
Then: $\_J_k$ is of full row rank; it has size $m_k\x n_k$ where $m_k\le n_k$; the sum of the differences $\sum_k (n_k-m_k)$ equals \DOF.
For each $k$, select $m_k$ columns of $\_J_k$ that form a nonsingular matrix $\_G_k$.
This can and must be done in such a way that the set of selected columns increases with $k$, so that each $\_G_k$ is a sub-matrix of the next.
For each of the $(n_k{-}m_k)$ {\em unselected} columns $j$ consider the item $\der{x}{j}{k+d_j}$.
The set of all these is a valid state vector $\_x_S$ since, briefly, non-singularity of $\_G_k$ ensures that at stage $k$, ``selected'' items $\der{x}{j}{k+d_j}$ belonging to selected columns can, by the Implicit Function Theorem, be found locally as functions of the unselected items. 

\smallskip
As said, \rf{livederiv} thus becomes a size-\DOF ODE system,
\begin{align}\label{eq:index1DAE1}
  \.1{\_x}_S = \_F(t,\_x_S).
\end{align}
This is locally equivalent to the size-$N_x$ DAE \rf{fider,livederiv} and hence to the original DAE.
Though ``index-1'' is the usual term used, the stronger property holds that
\bc \rf{fider,livederiv} form an {\em implicit ODE}, \ec
defined as an SA-friendly DAE whose offsets $c_i$ are all zero.

\subsection{Practical considerations}

At a DD-switch, the set \rf{fider} of differentiated equations does not change. 
Thus at the housekeeping level, a switch merely changes the set $S$ of index pairs $(j,l)$ that define the state vector.

We verified that the method works, by a proof-of-concept \matlab implementation. 
One example was the double pendulum (one pendulum-rod hung off another) in $x,y$ coordinates, where each rod independently has four DD-switching points in a full rotation, one in each of the four quadrants, giving $4\x4=16$ possible ``DD modes''.

It remains to be seen however how efficient one can make DD-switching for production code and for larger problems.
Finding the $\_G_k$ at a switch is nontrivial.
Ideally one wants each one to be maximally well-conditioned, which is expensive, so one seeks heuristic methods.
For this reason Scholz and Steinbrecher's simplified method \cite{scholz2013combined} is interesting.
It finds a state vector based on a highest-value transversal of the \sigmx; it is less general than full DDs but cheaper.
One might try it first, and if it gives ill-conditioned $\_G_k$ then use full DDs.

\smallskip
It seems natural to solve the original DAE numerically by giving formulation \rf{fider,livederiv} to a standard index-1 DAE solver.
However many models, especially mechanical ones, have many equations but few degrees of freedom, $N_x\gg\DOF$.
Then it makes sense to convert to the explicit form \rf{index1DAE1}. In many mechanical contexts (though not all) this ODE is typically non-stiff and thus amenable to solving by, say, an explicit Runge--Kutta method.
Working memory for sub-problems of size up to $n$ is needed by the root-finding that forms \rf{index1DAE1}, but is typically less than that needed by an implicit DAE code on a problem of size $N_x$.
For more discussion see \cite{Pryce2015c}.

\subsection{Example}

\begin{example}{Pendulum}
Let the original DAE be the simple pendulum in cartesian coordinates, shown with its \sigmx \rf{sigmx}, with relevant transversals marked.
Gravity $g$ and length $\ell$ are constants, and $x(t)$, $y(t)$ and $\lam(t)$ are state variables.
{\nc\s[1]{\mbox{\footnotesize\ensuremath{#1}}}
\rnc\o{^\circ}
\rnc\b{^\bullet}
\begin{align}\label{eq:pend0}
 & \begin{array}[m]{l}
    0 =    A = \.2x + x\lam, \\
    0 =    B = \.2y + y\lam - g, \\
    0 =    C = x^2 + y^2 - \ell^2,
  \end{array}
&& \Sigma =
\begin{blockarray}{rcccl}
    &x &y & \lam & \s{c_i}\\
 \begin{block}{r[ccc]l}
  A & 2\b&   & 0\o &\s0 \\
  B &   & 2\o & 0\b & \s0 \\
  C & 0\o & 0\b &   & \s2 \\
\end{block}
\s{d_j}          &\s2&\s2 & \s0 
\end{blockarray}
\end{align}
}\vspace{-2ex}

To describe DDs here, we modify the general $x_{jl}$ notation by renaming $x,\.1x,\.2x$ to $x_0,x_1,x_2$, and so on, to get 5 equations in 7 unknowns:
\begin{align}\label{eq:pend1}
 & \begin{array}[m]{l}
  \mc1c{\text{Augmented system}} \\
    0 =    A = \.2x + x\lam \\
    0 =    B = \.2y + y\lam - g \\
    0 =    C = x^2 + y^2 - \ell^2 \\
    0 = \.1C = 2(x\.1x + y\.1y) \\
    0 = \.2C = 2(x\.2x + \.1x^2 + y\.2y + \.1y^2) \\[1ex]
    \mc1l{\text{unknowns $x,\.1x,\.2x,  y,\.1y,\.2y, \lam$}}
  \end{array}
&&
  \begin{array}[m]{l}
  \mc1c{\text{After renaming}} \\
    0 =  A_0 = x_2 + x_0\lam_0 \\
    0 =  B_0 = y_2 + y_0\lam_0 - g \\
    0 =  C_0 = x_0^2 + y_0^2 - \ell^2 \\
    0 =  C_1 = 2(x_0x_1 + y_0y_1) \\
    0 =  C_2 = 2(x_0x_2 + x_1^2 + y_0y_2 + y_1^2) \\[1ex]
    \mc1l{\text{unknowns $x_0,x_1,x_2, y_0,y_1,y_2, \lam_0$}}
  \end{array}
\end{align}
One can choose any of $(x,\.1x),(y,\.1y),(x,\.1y),(y,\.1x)$ as state vector (one {\em must} choose one undifferentiated variable and one first derivative), but only the first two are ``convenient'' for AD, as the next section shows.

Suppose for example $\_x_S = (x,\.1x) = (x_0,x_1)$.
It is easily seen that provided $y\ne0$ one can find all the items as functions of these two, hence the pendulum DAE is equivalent to an ODE \rf{index1DAE1} in this $\_x_S$ when $y\ne0$.

The description of DDs given in \ssrf{DDs} has the advantage of combining {\em index reduction} and {\em order reduction} into one process.
For computer solution, it is probably easiest to work with the order 1 DAE formed by the $N_x=n+\sum_j d_j$ equations \rf{fider,livederiv}.
However ``by hand'', one can simplify by directly substituting the derivative relations into \rf{fider} where possible. E.g., when $\_x_S$ is $(x_0,x_1)$, one obtains
\begin{align*}
  0 &= A_0 = \.1x_1 + x_0\lam_0 \\
  0 &= B_0 = y_2 + y_0\lam_0 - g \\
  0 &= C_0 = x_0^2 + y_0^2 - \ell^2 \\
  0 &= C_1 = 2(x_0x_1 + y_0y_1) \\
  0 &= C_2 = 2(x_0\.1x_1 + x_1^2 + y_0y_2 + y_1^2) \\
  0 &= \.1x_0 - x_1
\end{align*}
The last equation, $\.1x_0 = x_1$, can not be ``substituted away''---in general, any equation \rf{livederiv} must stay if its $x_{jl}$ and $x_{j,l+1}$ are both state items, as this is how order reduction occurs.

\smallskip
In \cite{Matt93a}'s terminology a ``dummy derivative'' means a differentiated item that, in our terms, is a solved for item but is not a state variable or the derivative of one. In this example with this state vector, that makes $y_1$ and $y_2$ the DDs.
\end{example}

\subsection{AD for dummy derivatives}\label{ss:howADhelp}

How can an AD tool help automate numerical solution by DDs, as described above?

First, it is helpful if the tool supports $\d/\d t$ as a first-class operator, of equal status with $+,\x,\sin()$, etc., so that it can understand a representation of a DAE in the general form \rf{maineq}.
This is not essential. Tools such as ADOL-C and dcc/dco \cite{ADOLC,Naumann2012TAo} do not have this feature, but can handle arbitrary expressions containing derivatives by renaming the latter as algebraic items and stating their differential relations separately.
This is like the method in \ssrf{DDs}, where derivatives are renamed as algebraic in \rf{xjder} and some differential relations between them stated in \rf{livederiv}.

Our solver \daets uses Ole Stauning's AD package \fadbadpp \cite{FADBAD++}, written in \cpp.
It did not originally include $\d/\d t$ but at our request in 2002, Stauning included the operator \@Diff@ such that  \@Diff($\cdot,q$)@ means $\d^q/\d t^q$.
For instance, straightforward code for the pendulum, as in the \daets user guide, is shown in \fgref{pendprog}.
\begin{figure}[ht]
\begin{lstlisting}
template <typename T>                /*\label{ln:a}*/
void fcn(T t, const T *z, T *f, void *param) {
 // z[0], z[1], z[2] are $x$, $y$, $\lambda$.
 const double G = 9.81, L = 10.0;
 f[0] = Diff(z[0],2) + z[0]*z[2];
 f[1] = Diff(z[1],2) + z[1]*z[2] - G;
 f[2] = sqr(z[0]) + sqr(z[1]) - sqr(L);
}                                    /*\label{ln:b}*/
\end{lstlisting}
  \caption{Code for simple pendulum problem. \label{fg:pendprog}}
\end{figure}

\smallskip
More important, for DDs and other index reduction methods, an AD tool must be able to differentiate the $f_i$ selectively. For instance in the pendulum, $A$ and $B$ are to be left alone, and $C$ differentiated twice.

At first sight this seems to require a tool based on source code transformation, which could generate code symbolically for the last two equations in \rf{pend1}, for instance.
But this is not the case---the key is not to treat different derivatives of a given variable in isolation, but store them together as a truncated power series.
For instance in the pendulum, the unknowns form three objects
\[\begin{array}{rll}
     x &= (x_0,x_1,x_2) &\text{order 2 power series}, \\
     y &= (y_0,y_1,y_2) &\text{order 2 power series}, \\
  \lam &= (\lam_0) &\text{order 0 power series}.
\end{array}\]
To keep the algebra simple below, these are Taylor coefficients, not derivatives, thus $x_k$ here is the $x_k$ in \rf{pend1} divided by $k!$, and so on.

Taylor series AD by overloading, provided by many AD tools, now gives the needed values. For instance evaluating $C=x^2+y^2-\ell^2$ proceeds via these intermediate steps:
\[\begin{array}{rrll}\hline
  \mc2l{\text{input}} \\
  x &&= (x_0,x_1,x_2) \\
  y &&= (y_0,y_1,y_2) \\\hline
  \mc2l{\text{compute}} \\
  v_1&= x^2 &= (x_0^2,\quad 2x_0x_1,\quad 2x_0x_2+x_1^2) \\
  v_2&= y^2 &= (y_0^2,\quad 2y_0y_1,\quad 2y_0y_2+y_1^2) \\
  v_3&= v_1+v_2 &= (x_0^2+y_0^2,\quad 2(x_0x_1+y_0y_1),\quad 2(x_0x_2+y_0y_2) + x_1^2+y_1^2) \\
  C  &= v_3 - \text{const}(\ell^2) &=
              (x_0^2+y_0^2-\ell^2,\quad 2(x_0x_1+y_0y_1),\quad 2(x_0x_2+y_0y_2) + x_1^2+y_1^2) 
\end{array}\]
returning an order 2 power series object $C$ holding the needed coefficients $(C_0,C_1,C_2)$, that is $(C,\.1C,\frac12\.2C)$ in terms of derivatives.

Evaluating $A = \.2x + x\lam$ and $B = \.2y + y\lam - g$ is similar. $\@Diff@(\cdot,2)$ converts, e.g., the order 2 series $x = (x_0,x_1,x_2)$ to the order 0 series $(2x_2)$.
Thus $A,B$ are returned as order 0 series $A=(A_0)=(2x_2+x_0\lam_0)$, $B=(B_0)=(2y_2+y_0\lam_0-g)$.

\smallskip
The above method gives an explicit evaluation of the $N_f$ functions \rf{fider} at the $N_x$ arguments \rf{xjder}.

In the context of DDs and reducing the DAE to an explicit ODE, state item values, say values $\_x_S=(x_0,x_1)$ are given; the 5 items $\_x_F=(x_2,y_0,y_1,y_2,\lam_0)$ are trial values that produce 5 residual values $\_r=(A_0,B_0,C_0,C_1,C_2)$. By root-finding using suitable Jacobians (found by methods not described here) we find $\_x_F$ that makes $\_r=\_0$, thus solving for $\_x_F$ as a function of $\_x_S$. Extract $x_2$ from $\_x_F$ to form $(x_1,x_2)$, which is $\.1{\_x}_S$.
This implements the function $\_F$ in \rf{index1DAE1}.

To make this work, the state items must comprise a contiguous set of derivatives of each variable, with no gaps.
(Hence, see below \rf{pend1}, $(x,\.1y)$ and $(y,\.1x)$ are not useful state vectors for the pendulum.)
That is, $S$ must have the form $\set{(j,l)}{$0\le l<\delta_j$, $j=1,\ldots,n$}$, where $\bdelta=(\delta_1,\ldots,\delta_n)$ is an integer {\em DD-spec vector} with $0\le\delta_j\le d_j$ and $\sum_j\delta_j=\DOF$, which uniquely specifies the DD scheme currently in use. DD switching can be based on changing this $\bdelta$, and following through the consequences for various associated index sets and Jacobian-related matrices.

Experiments by Nedialkov confirm this is an effective and flexible way to implement DDs, including the root-finding that produces \rf{index1DAE1}, using methods already provided by the \daets classes.

\section{The Lagrangian}\label{sc:lagrangian}

\subsection{Lagrangian mechanics theory}\label{ss:lagtheory}

For mechanical systems, such as in robotics, equations of motion can often be conveniently derived from the system's Lagrangian function $L$.
It is assumed there are conservative (energy preserving) forces such that one can define a potential energy $V$ depending only on system position.
Then $L=T-V$, where $T$ is the system's total kinetic energy.
Let the configuration at any time be described by generalised coordinates $\_q = (q_1,\ldots,q_{n_q})$ such that $T$ is a function of $\.1{\_q}$ and possibly $\_q$, and $V$ is a function of $\_q$ only. There may (depending on the coordinate system used) be $n_c$ scalar constraints that are holonomic, i.e.\ functions of positions and possibly time but not of velocities, namely $C_j(t,\_q)=0$.
Then the variational principle of stationary action gives the $(n_q{+}n_c)$ equations of motion:
\begin{align}
  \frac{\d}{\d t} \dbd{L}{\.1q_i} - \dbd{L}{q_i}
   + \sum_{j=1}^{n_c} \lam_j \dbd{C_j}{q_i} &=0,
  \quad i=1,\ldots,n_q, \label{eq:lagrange1}\\
  C_j(t,\_q)&=0, \quad j=1,\ldots,n_c, \label{eq:lagrange2}
\end{align}
where the $\lam_j$ are Lagrange multipliers for the constraints. If $n_c>0$, i.e.\ constraints are present, \rf{lagrange1,lagrange2} is termed a Lagrangian system {\em of the first kind}. It is a DAE system, of index 3 in the classical sense or index $2$ as defined in \rf{indexdef}, since two $t$-differentiations of each $C_j$ are needed.
If the coordinates are chosen so that $n_c=0$, it is {\em of the second kind} and is an ODE system.

For a system subject to external forces, the zero right hand sides of \rf{lagrange1} are replaced by $u_i(t)$, $i=1,\ldots,n_q$, which are {\em generalised external force components}.
One can model certain kinds of internal dissipative (not energy conserving) velocity-dependent force by a so called {\em Rayleigh dissipative term}.

For instance, for free motion of the simple pendulum, taking $\_q$ to be $(x,y)$, the coordinates of the pendulum bob with $y$ measured downward, gives
\begin{align}\label{eq:pendl}
  L = \frac12 m(\.1x^2+\.1y^2) + mgy, 
\end{align}
where $m$ is the bob's mass, and with one constraint $x^2+y^2-\ell^2=0$. 
Then \rf{lagrange1,lagrange2} lead to equations \rf{pend0}.
But taking  $\_q$ to be $(\theta)$, the angle of the pendulum from the downward vertical, gives $L = \frac12 m(l\.1\theta)^2 + mgl\cos\theta$, with no constraints. 
Then \rf{lagrange1,lagrange2} lead to the ODE $\.2\theta = -(g/l)\sin\theta$, which is equivalent to \rf{pend0}.

\medskip
Our method applies to solution by the \cpp code \daets.
We describe how it looks to the user, not the inner workings.
The user describes a DAE system to \daets using a template type \li{T}, as \fgref{pendprog} in \ssrf{howADhelp} shows. 
\daets instantiates \li{T} with several concrete types during execution to perform \SA and other tasks---in particular Taylor coefficient generation during numerical solution using \fadbadpp's Taylor type.

The Lagrangian facility overlays the \fadbadpp backward (reverse mode) differentiation type \li{B} on top of \li{T} to create the type \li{B<T>} which is then suitably manipulated.
It avoids the seemingly needed two ``phases'': the equations of motion are created, and readied for numerical solution, in one operation that happens at runtime.

\begin{figure}
\begin{lstlisting}
template <typename T>
void fcn( T t, const T *z, T *f, void *param ) {
  vector< B<T> > q(SIZEOFQ), qp(SIZEOFQ), C(SIZEOFC);               /*\label{pend:decl}*/
  init_q_qp(z, q, qp);                        /*\label{pend:init}*/
  double *p = (double *)param;                         /*\label{pend:p}*/
  double  m = p[0], g = p[1], l = p[2];                /*\label{pend:r1}*/
  B<T> x = q[0], y = q[1], xp = qp[0], yp = qp[1];     /*\label{pend:r2}*/
  B<T> L = 0.5 * m * (sqr(xp) + sqr(yp)) + m * g * y;  /*\label{pend:L}*/
  C[0] = sqr(x) + sqr(y) - sqr(l);                     /*\label{pend:C}*/
  setupEquations(L, z, q, qp, C, f);        /*\label{pend:set}*/
}/*\label{pend:e}*/
\end{lstlisting}
\caption{Code to describe pendulum in Lagrangian form.\label{fg:lagpendcode}}
\end{figure}
The pendulum in $x,y$ form may be coded as the \li{fcn} function in \fgref{lagpendcode} that defines a DAE in the form accepted by \daets.
Here \li{SIZEOFQ} and \li{SIZEOFC} are macros, set to $2$ and $1$.
The physical parameters $m$, $g$ and $\ell$ are passed in the array \li{param}. \daets specifies this as having type \li{void}, which leaves the user free to use some other type than \li{double} if needed.
(Since $m$ is a factor all through $L$, it is superfluous and could have been omitted.)

Lines \ref{pend:p}--\ref{pend:r2} name things for readability; lines \ref{pend:L}--\ref{pend:C} then define the problem.
Lines \ref{pend:decl}--\ref{pend:init} and \ref{pend:set} are boilerplate  and do not change from one problem to the next. The call \li{init_q_qp} ``connects'' the input variables of \rf{pend0} in \li{z} to the \li{q} and \li{qp} variables and prepares them for the partial derivatives computation.
The \li{setupEquation} call sets up the AD for converting the form \rf{pendl} to \rf{pend0}.
When \li{setupEquation} returns, \li{f[0]}, \li{f[1]}, and \li{f[2]} contain $A$, $B$ and $C$ in \rf{pend0}, (modulo a factor $m$).

\subsection{Examples}\label{ss:lagexamples}

We have applied the \daets Lagrangian facility to various systems, including the following examples. Performance tests are on a 2015 MacBook Pro laptop with a 4-core 2.2 GHz Intel processor running Mac OS X 10.11.6.
The \cpp compiler is  GCC 6.3.0.
All numerics are in \cpp \li{double}.

Nedialkov's group at McMaster University have recently improved substantially the performance of the AD in the Taylor series method of \fadbadpp by applying common subexpression elimination techniques.
They have also implemented an efficient algorithm for computing the System Jacobians,  propagating compressed gradients, and have incorporated sparse linear algebra in \daets. 
Speedups of $>10$ compared to using the original \fadbadpp are common. 
As a result code created by the Lagrangian facility is surprisingly fast.

\begin{wrapfigure}[14]{r}{.35\TW}
  \includegraphics[width=.35\TW]{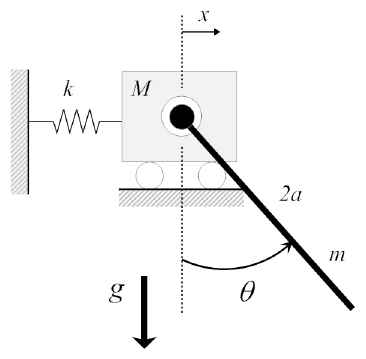}
  \caption{Spring-Mass-Pendulum with one rod.}
\end{wrapfigure}
Because of the perceived difficulty of solving DAEs, generalised coordinates are often chosen to eliminate the constraints and give a Lagrangian of the second kind.
For instance a rigid body's 3D position can be described by 3 coordinates of the position of its \CM and 3 of its angular position relative to this. 
However the mathematical formulation is often simpler in cartesian coordinates. 
One plus of using a code for high-index DAEs such as \daets is that it handles resulting ``first kind'' systems easily.
Further, since \daets does not set up a local coordinate system for numerical solution as the DDs method does, it does not suffer the performance penalty of DD-switching.

\begin{example}{Spring-Mass-Pendulum}
\def\r#1dot{\.1{\_r}_{#1}}
\def\x#1dot{\.1x_{#1}}
\def\y#1dot{\.1y_{#1}}
\nc\m{\textup{\;m}}
\nc\s{\textup{\;s}}
\nc\kg{\textup{\;Kg}}
This 2D model is taken from an article on the Acumen mechanics modeling system by Zhu, Taha {\it et al.} \cite{ZhuTaha2010}.

We have extended their model to a chain of any number $n$ of rods.
Namely, a horizontally sliding point-mass $M$ is connected by a spring of stiffness $k$ to a fixed point at the same level.
From $M$ hangs a chain of $n$ uniform rods, with frictionless joints between the end of one and the start of the next.
Purely to simplify the code, they all have the same mass $m$ and length $l=2a$.
We assume the setup is constructed so that all components can slide or rotate freely without colliding.

For $n\ge2$ (possibly even for $n=1$) the motion can be chaotic.
The figure (taken from \cite{ZhuTaha2010}) shows the case $n=1$.
As the figure indicates, \cite{ZhuTaha2010} takes $\_q=(x,\theta)$ as coordinates, leading to a Lagrangian of the second kind, $L=T-V$ where:
\begin{align}\label{eq:zhutaha1}
   T &= \tfrac12(M+m)\.1x^2 + ma\.1x\.1\theta\cos\theta + \tfrac23 ma^2 \.1\theta^2,
  &V &= \tfrac12 kx^2 +mga(1-\cos\theta).
\end{align}
Here the rotational kinetic energy term $\tfrac23 ma^2 \.1\theta^2$ uses the moment of inertia $I=\frac43 ma^2$ of a uniform rod about its \CM.

For the general $n$-rod model we use, instead, cartesian coordinates $\_q = (x_0,x_1,\ldots,x_n;\; y_1,\ldots,y_n)$. Here $\_r_0 = (x_0,y_0)$, with $y_0$ constant equal to $0$, is the position of $M$ and the start of rod $1$, and $\_r_i = (x_i,y_i)$ for $i=1,\ldots,n$ is the position of the joint between the end of rod $i$ and (for $i<n$) the start of rod $i+1$.
We avoid moments of inertia by using the following, where $\cdot$ denotes the dot product of vectors.
\begin{lemma}
  If the ends of a uniform rod of mass $m$ have position vectors $\_r_0$ and $\_r_1$, depending on $t$, then its kinetic energy at any instant is
\[
  \textup{KE} = \tfrac16 m (\r0dot\cdot\r0dot + \r0dot\cdot\r1dot + \r1dot\cdot\r1dot).
\]
\end{lemma}
\begin{proof}
  We can parameterise position along the rod as $\_r=(1-s)\_r_0 + s\_r_1$, for $0\le s\le1$.
  Since the rod has total mass $m$, an element from $s$ to $s+\d s$ has mass $m\d s$.
  The velocity of this element is $\.1{\_r}=(1-s)\.1{\_r}_0 + s\.1{\_r}_1$ so its kinetic energy is
  \begin{align*}
  \tfrac12 m (\.1{\_r}\cdot\.1{\_r}) \d s 
  = \tfrac12 m \left((1-s)^2\;\.1{\_r}_0\cdot\.1{\_r}_0
    + 2(1-s)s\;\.1{\_r}_0\cdot\.1{\_r}_1
    + s^2\;\.1{\_r}_1\cdot\.1{\_r}_1 \right)\d s.
  \end{align*}
  Integrating this from 0 to 1 gives the result.
\end{proof}
The potential energy of the rods comes from considering the mass of rod $i$ to be at its \CM at height $\frac12(y_{i-1}+y_i)$; there is a contribution of $\frac12 k x_0^2$ from the spring and none from mass $M$.
This leads to the Lagrangian $L=T-V$, and constraints $C_i$, where
\begin{align*}
  T &= \tfrac12 M \x0dot^2 + \tfrac16 m \sum_{i=1}^n 
    \left(\r{i-1}dot\cdot\r{i-1}dot + \r{i-1}dot\cdot\r{i}dot + \r{i}dot\cdot\r{i}dot \right) \\
    &= \tfrac12 M \x0dot^2 + \tfrac16 m \sum_{i=1}^n 
    \left((\x{i-1}dot^2 + \y{i-1}dot^2) + (\x{i-1}dot\x{i}dot + \y{i-1}dot\y{i}dot)
    +  (\x{i}dot^2 + \y{i}dot^2)\right) \\
  V &= \tfrac12 k x_0^2 + mg \sum_{i=1}^n \tfrac12(y_{i-1}+y_i)
     = \frac12 k x_0^2 + mg \Bigl(\tfrac12 y_n + \sum_{i=1}^{n-1} y_i\Bigr), \\
  0=C_i &= (x_i-x_{i-1})^2 + (y_i-y_{i-1})^2 - \ell^2, \quad
   (i=1,\ldots,n).
\end{align*}
The code in \fgref{ex1code}, which replaces lines 
\ref{pend:r2}--\ref{pend:C} in the \li{fcn} of \fgref{lagpendcode}, implements the above formulas. 
Here \li{n}, the number $n$ of rods, is read in as one of the physical parameters. \li{SIZEOFC} also equals $n$.
The arrays \li{q} and \li{qp} holding $\_q$ and $\.1{\_q}$ have length $2n+1$.

\begin{figure}
\begin{lstlisting}
B<T> x0 = q[0], *x = q+1, *y = x+n, x0p = qp[0], *xp = qp+1, *yp = xp+n;
  B<T> L; 
  {                                              /*\label{spring:b1}*/
    B<T> KEsum = sqr(x0p) + sqr(xp[n-1]) + sqr(yp[n-1]) + x0p*xp[0];
    for (int i=0; i < n-1; i++)
      KEsum += 2*( sqr(xp[i]) + sqr(yp[i]) ) + xp[i]*xp[i+1] + yp[i]*yp[i+1];
    B<T> KE = 0.5*M*sqr(x0p) + m/6*KEsum;
   
    B<T> PEsum = 0.5*y[n-1];
    for (int i=0; i<n-1; i++) PEsum += y[i];
    B<T> PE = 0.5*k*sqr(x0) - m*g*PEsum;   // /*$-$ as $y$ goes downward */
   
    L = KE - PE;
  }                                            /*\label{spring:b2}*/
  C[0] = sqr(x[0]-x0) + sqr(y[0]) - sqr(l);
  for (int i=1; i<SIZEOFC; i++)
    C[i] = sqr(x[i]-x[i-1]) + sqr(y[i]-y[i-1]) - sqr(l);  
\end{lstlisting}
\caption{Code for spring-mass-multi-pendulum system.
\label{fg:ex1code}}
\end{figure}

Listing line 1 uses \li{C} syntax\footnote{For instance \li{q+1} references a sub-array of \li{q} starting at \li{q[1]}.}
 to split \li{q} into a scalar holding $x_0$, and two size-$n$ arrays holding $x_1,\ldots,x_n$, and $y_1,\ldots,y_n$; similarly \li{qp}.
The variable \li{KEsum} accumulates 
$\dot{x}_0^2 + \dot{x}_n^2+\dot{y}_n^2+ \dot{x}_0\dot{x}_1+ \sum_{i=1}^{n-1}\bigl[2(\dot{x}_i^2+\dot{y}_i^2)+\dot{x}_i\dot{x}_{i+1}+ \dot{y}_i\dot{y}_{i+1}\bigr]$, which is equivalent to the sum in $T$, and similarly \li{PEsum}.

In the computation of \li{L}, the temporary dependent variables \li{KEsum}, \li{PEsum}, \li{KE}, and \li{PE}, are local in  the block between lines \ref{spring:b1} and \ref{spring:b2}; \fadbadpp requires that in the reverse mode either all  intermediate dependent variables are differentiated or go out of scope, which is the case here. 

In our tests the physical parameters of the original model in \cite{ZhuTaha2010} were used, namely, assuming SI units, $g=9.8\m\s^{-2}$, $l=2a=2\m$, $M=5\kg$, $m=2\kg$, $k=10\kg\s^{-2}$.

The chosen initial conditions (ICs) are that the system is at rest with mass $M$ at $x=4$, and the rods stretched horizontally to the left.
(Thus the spring is pushing against the row of rods; animations show it ``folds up'' rods 1 and 2 as they start to fall.)

\smallskip
To confirm that we are modeling the same system as in \cite{ZhuTaha2010}, the equations of motion derived from the Lagrangian \rf{zhutaha1} given in \cite{ZhuTaha2010} were coded in \matlab and integrated by \li{ode45}.
The results were compared with those of the \daets version for the case $n=1$.
The latter was coded to output $q$ and $\.1q$ at each of its time points $t_i$. This data was mapped to the $t_i$ chosen by the \matlab version by Hermite cubic interpolation between adjacent $t_i$ of \daets.
\fgref{ex1errplot} shows that over $t=[0,40]$, the differences ({\li{ode45} solution at tolerance $10^{-12}$}) $-$ ({\daets solution at tolerance $10^{-8}$)} are of order $10^{-6}$. This gives confidence that the programs are solving the same physical model.
\begin{figure}\centering
  \includegraphics[width=.6\TW]{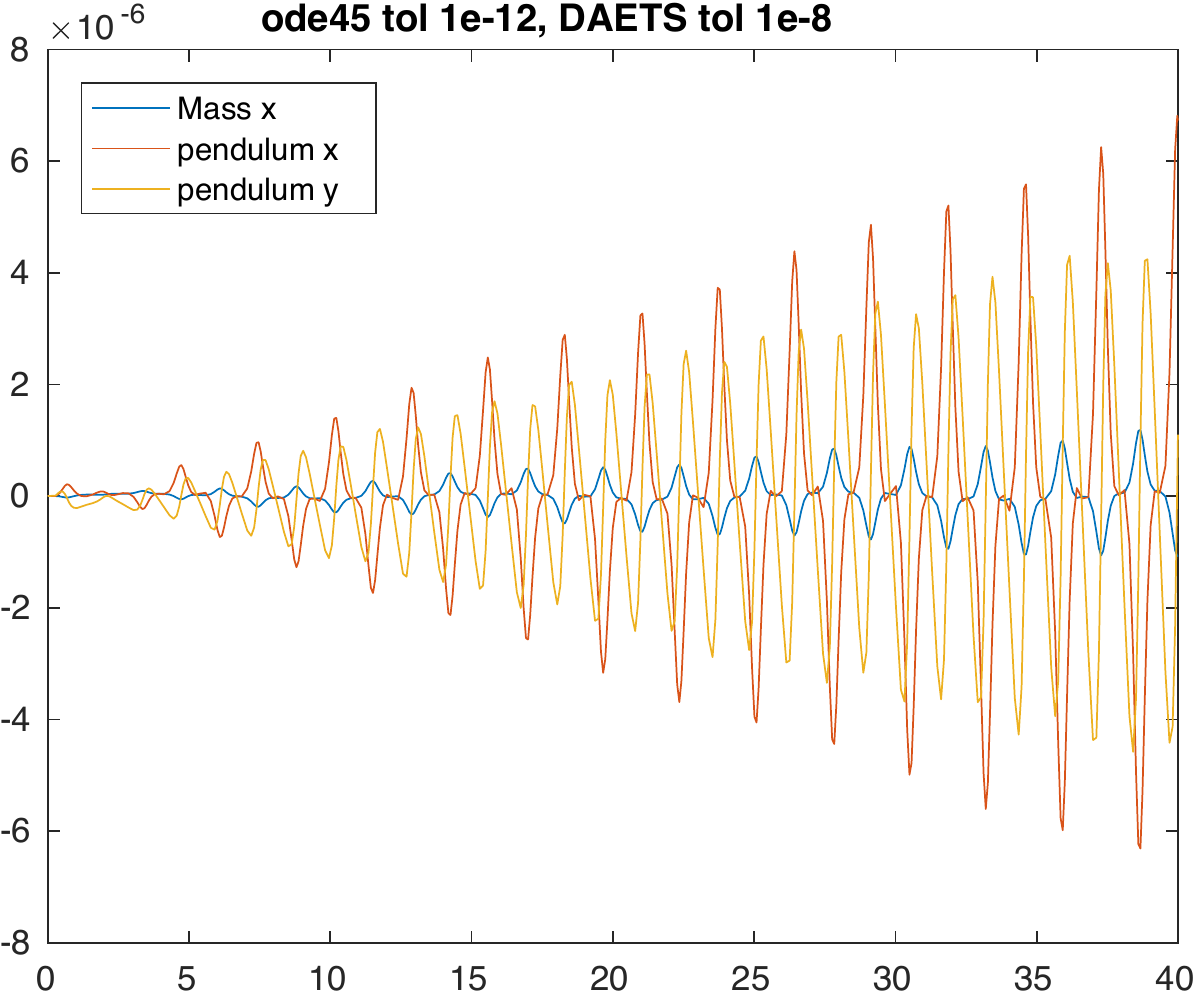}
  \caption{Spring-mass-pendulum with $n=1$.
  It shows, for the $x$ coordinate of the sliding mass and the $x$ and $y$ coordinates of the pendulum, the difference between the solution by our model and that by \protect\rf{zhutaha1}, over $0\le t\le40$.\label{fg:ex1errplot}}
\end{figure}

\begin{table}\small
\begin{tabular}{c|ccccrrccccc}
 \@tol@ &$n=$1 &2    &4    &6    &\multicolumn{1}{c}{8}     &\multicolumn{1}{c}{10}    &12    &14    &16    &18    &20 \\\hline
 1e--4  &0.19 &0.44 &1.06 &2.98 & 6.12 & 8.94 &12.10 &19.75 &29.72 &39.83 &47.96\\
 1e--6  &0.23 &0.56 &1.45 &4.49 & 8.11 &11.03 &17.70 &25.28 &37.37 &53.56 &65.68\\
 1e--8  &0.34 &0.68 &1.83 &5.76 &10.84 &14.94 &24.13 &35.42 &51.12 &73.42 &86.13\\
 1e--10 &0.40 &0.91 &2.63 &7.41 &14.92 &19.02 &33.29 &44.75 &67.31 &98.79 &121.5\\
 1e--12 &0.55 &1.35 &3.59 &9.99 &18.40 &26.89 &45.90 &62.79 &96.56 &124.5 &155.9\\[1ex]
\end{tabular}
\caption{Time (seconds) to integrate to $t=100$ for various numbers $n$ of rods, and tolerances \protect\li{tol}.\label{tb:ex1timings}}
\end{table}

\smallskip
For timing tests, the system was integrated by \daets over $0\le t\le100$ for various numbers $n$ of rods and (mixed relative-absolute) tolerances \li{tol}.
Case $n=1$ is the model in \cite{ZhuTaha2010}.
The Taylor series order was set to 15, which works well for these problems at this range of accuracies.
\tbref{ex1timings} shows the times taken.

\daets has a ``maximum step size'' feature but this was not used so it chose the step sizes $h$ freely. 
For the ``hardest'' problem $n=20$ at tolerance $10^{-12}$ they ranged from $h=0.00061$ to $h=0.093$.
For the ``easiest'', $n=1$ at tolerance $10^{-4}$, they ranged from $h=0.065$ to $h=0.421$.
\end{example}

\begin{example}{Controlled simple pendulum}

We show one can solve a prescribed-trajectory control problem for a Lagrangian-described system.
Namely, for the simple pendulum we introduce a horizontal external force on the bob, modeled as a system input $u = u(t)$ such that the equation $\ddot x + \lam x = 0$ becomes
\begin{align}\label{eq:pendsysin}
  \ddot x + \lam x - u = 0.
\end{align}
The aim is to find $u(t)$ (plus suitable consistent ICs) so that the $x$ position performs simple harmonic motion 
$x(t) = a\sin(\omega t)$
exactly, where the constants $a$ and  $\omega$ are a given amplitude and frequency, respectively.

Comparing the pendulum as initial-value problem in \fgref{lagpendcode} and as control problem in\fgref{lagpendcode2} shows the {\em implementation} changes little. One passes $a$ and $\omega$ as extra parameters that become \li{a} and \li{w}.
After the \li{setupEquations} line, the first equation \li{f[0]} is modified in line \ref{pendf:x}, and a new fourth equation \li{f[3]} is at line \ref{pendf:u} (in which the \li{x} at line \ref{pend:r2} in \fgref{lagpendcode} cannot be used as it has the wrong type, \li{B<T>} instead of \li{T}).

But the revision has changed the DAE's mathematical nature greatly.
Now with $4$ variables and equations, it is shown below with its signature matrix $\Sigma$ (a blank means $\ninf$, and the unique transversal is marked by $^\circ$).
{\nc\s[1]{\mbox{\footnotesize\ensuremath{#1}}}
\rnc\o{^\circ}
\begin{align}\label{eq:pendcontrol}
\renewcommand{\arraystretch}{1.3}
&\begin{array}{l}
 0 = A = x'' +x \lam - u \\
 0 = B = y''+y\lam-g \\
 0 = C = x^2+y^2-\ell^2 \\
 0 = D = x - a\sin(\omega t) \\
 ~
\end{array}
&& \Sigma =
\begin{blockarray}{rccccl}
      &x  &y &\lam& u & \s{c_i}\\
     \begin{block}{r[cccc]l}
      A & 2 &   & 0 &0\o& \s0 \\
      B &   & 2 &0\o&   & \s0 \\
      C & 0 &0\o&   &   & \s2 \\
      D &0\o&   &   &   & \s2 \\
    \end{block}
\s{d_j} &\s2&\s2&\s0&\s0&
\end{blockarray}
\end{align}
}
While \rf{pend0} has $2$ degrees of freedom, \rf{pendcontrol} has none---specifying the desired $x(t)$ determines the system input $u(t)$, as well as $y$ and $\lambda$, uniquely.
\begin{figure}
\begin{lstlisting}
template <typename T>
void fcn( T t, const T *z, T *f, void *param ) {
    vector< B<T> > q(SIZEOFQ), qp(SIZEOFQ), C(SIZEOFC);   /*\label{pendf:decl}*/
  init_q_qp(z, q, qp);                    /*\label{pendf:init}*/
  double *p = (double *)param;                /*\label{pendf:p}*/
  double  m = p[0], g = p[1], l = p[2],       /*\label{pendf:r1}*/
    a = p[3], w = p[4];                       /*\label{pendf:aw}*/
  B<T> x = q[0], y = q[1], xp = qp[0], yp = qp[1];     /*\label{pendf:r2}*/
  B<T> L = 0.5 * m * (sqr(xp) + sqr(yp)) + m * g * y;  /*\label{pendf:L}*/
  C[0] = sqr(x) + sqr(y) - sqr(l);                     /*\label{pendf:C}*/
  setupEquations(L, z, q, qp, C, f);       /*\label{pendf:set}*/
  f[0]-= z[3];                                       /*\label{pendf:x}*/
  f[3] = z[0] - a*sin(w*t);                          /*\label{pendf:u}*/
}
\end{lstlisting}
	\caption{The \li{fcn} for the controlled pendulum problem: new lines \ref{pendf:x} and \ref{pendf:u} are inserted.\label{fg:lagpendcode2}}
\end{figure}

With the physical parameters $g=9.8$ and $\ell=10$, the problem was solved by \daets with $\omega$ equal to the pendulum's natural frequency $\sqrt{g/l}$ of small oscillations, and for various $a$; and again with $\omega$ changed by $20\%$, for the same $a$ values.
Some examples of resulting $u$'s are plotted over several cycles in \fgref{pendcontrol}.
As expected, $u$ is very small when $\omega$ is the natural frequency and $a$ is small. It becomes large as $a$ approaches $\ell$, or as the frequency moves away from the natural one.
\daets took less than 0.1 seconds for each of the runs.
\begin{figure}\centering
  \includegraphics[width=1\linewidth]{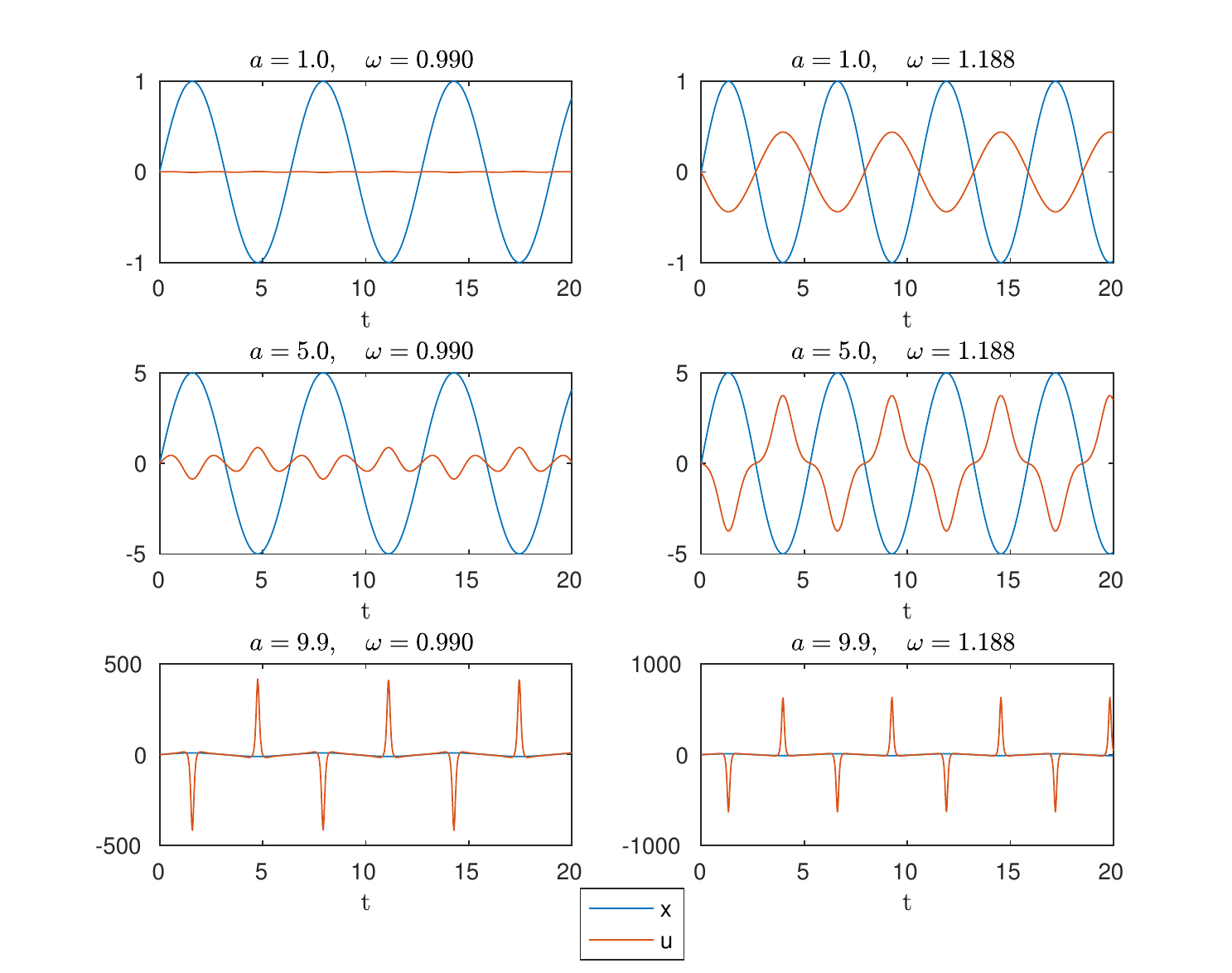}
  \caption{Solution by \daets of system input $u(t)$ for controlled pendulum with $g=9.8$, $l=10$. Required response $x = a\sin(\omega t)$. For $\omega$ equal to natural angular frequency (left column) and $20\%$ larger (right column), and three $a$ values.\label{fg:pendcontrol}}
\end{figure}
\end{example}

\begin{example}{DETEST Non-stiff Problem C5}\label{ex:DETEST}
\nc\brho{\bm{\varrho}}

This problem from the non-stiff part of the DETEST testing package for ODE solvers \cite{EnrightPryce1987}, and originally\footnote{\cite{EnrightPryce1987} cites nonexistent reference ``11'' which should be ``10'' and is the Zonneveld work.} from Zonneveld \cite{zonneveld1970automatic}, is titled ``Five Body Problem: Motion of five outer planets about the Sun''.
It is a order 2 ODE of size 15 (so size 30 when reduced to order 1), the variables being the positions of Jupiter, Saturn, Uranus, Neptune and Pluto relative to the Sun, in $x,y,z$ coordinates such that the ecliptic plane, in which the orbits approximately lie, is not close to any of the three axes.

To set up the Lagrangian formulation, $\_q$ comprising the 5 relative positions $(\brho_1(t),\ldots,\brho_5(t))$ (each $\brho$ being a 3-vector) is converted to 6 positions $(\_r_0(t),\ldots,\_r_5(t))$ of Sun and planets relative to their common \CM, which may be considered to be at rest in a Newtonian absolute frame.
Namely let $m_0$ be the mass of the Sun and $m_1,\ldots,m_5$ the masses of the planets and subtract
\[ \_r_c = \frac{m_0\_0 + (m_1\brho_1 + \cdots + m_5\brho_5)}{m_0\ + (m_1\ +\ \cdots\ +\ m_5)} \]
from each component of $(\_0,\brho_1,\ldots,\brho_5)$ to get $(\_r_0,\ldots,\_r_5)$.
Then
\begin{align}
  T &= \tfrac12\sum_{i=0}^5 m_i |\.1{\_r}_i|^2,
  &V &= -\sum_{i=0}^5\sum_{j=i+1}^5 \frac{Gm_im_j}{|\_r_i - \_r_j|},
  &L &= T-V,
\end{align}
where $G$ is the gravitational constant.
The code, shown in Appendix \ref{sc:nbody3Drelative}, was made particularly compact using a \cpp 3-vector class from \cite{vec3D2015}.

In the DETEST model the time unit (TU) is 100 days. Distance is measured in astronomical units (AU), where 1 AU is the mean radius of the earth's orbit.
The task is to integrate from given initial values up to $t=20$ TU; at tolerance $10^{-13}$ we get agreement with DETEST's reference solution to around 12 decimal places.

To see how fast the solution is, the problem was integrated to $t=200,000$ TU (about 55,000 earth years), with Taylor order 15, at tolerances $10^{-13}$ and $10^{-14}$. The two sets of results agree to five decimal places, and some \daets integration statistics are
\[
\begin{tabular}{c|cccccc}
 \mc5c{Integration of Sun and 5 planets to $t=200,000$ TU} \\\hline
 \@tol@   &CPU secs &no. of steps &smallest step & largest step \\\hline
 \@1e-13@ &31.2     &58028        &1.996         &5.654 \\
 \@1e-14@ &36.2     &67035        &1.725         &5.370
\end{tabular}
\]

It is known that Pluto is locked in a, currently stable, 3:2 resonance with Neptune.
This was easy to verify over short periods from our results.
For the subtleties of solar system behaviour see \cite{Hayes_surfingon2008} and references therein.
This article cites evidence that over very long times the system switches between regular and chaotic behaviour in an irregular way that depends critically on ICs.
Hence numerical results showing linear (regular) divergence of neighbouring solutions up to some large time $T$---rather than exponential (chaotic) divergence---are no evidence that such behaviour will continue up to, say, time $2\,T$.

What about the given ICs?
We integrated the problem at two tolerances \li{1e-13} and \li{1e-15}, recording the solutions at successive powers of ten up to $10^8$ TU ($<5$ hours CPU time for each to reach $10^8$) and computing the relative error in the two norm  at these times.
The results, see \fgref{planetdiverge}, indicate non-chaotic behaviour up to that point.
\begin{figure}\centering
  \includegraphics[width=0.6\TW]{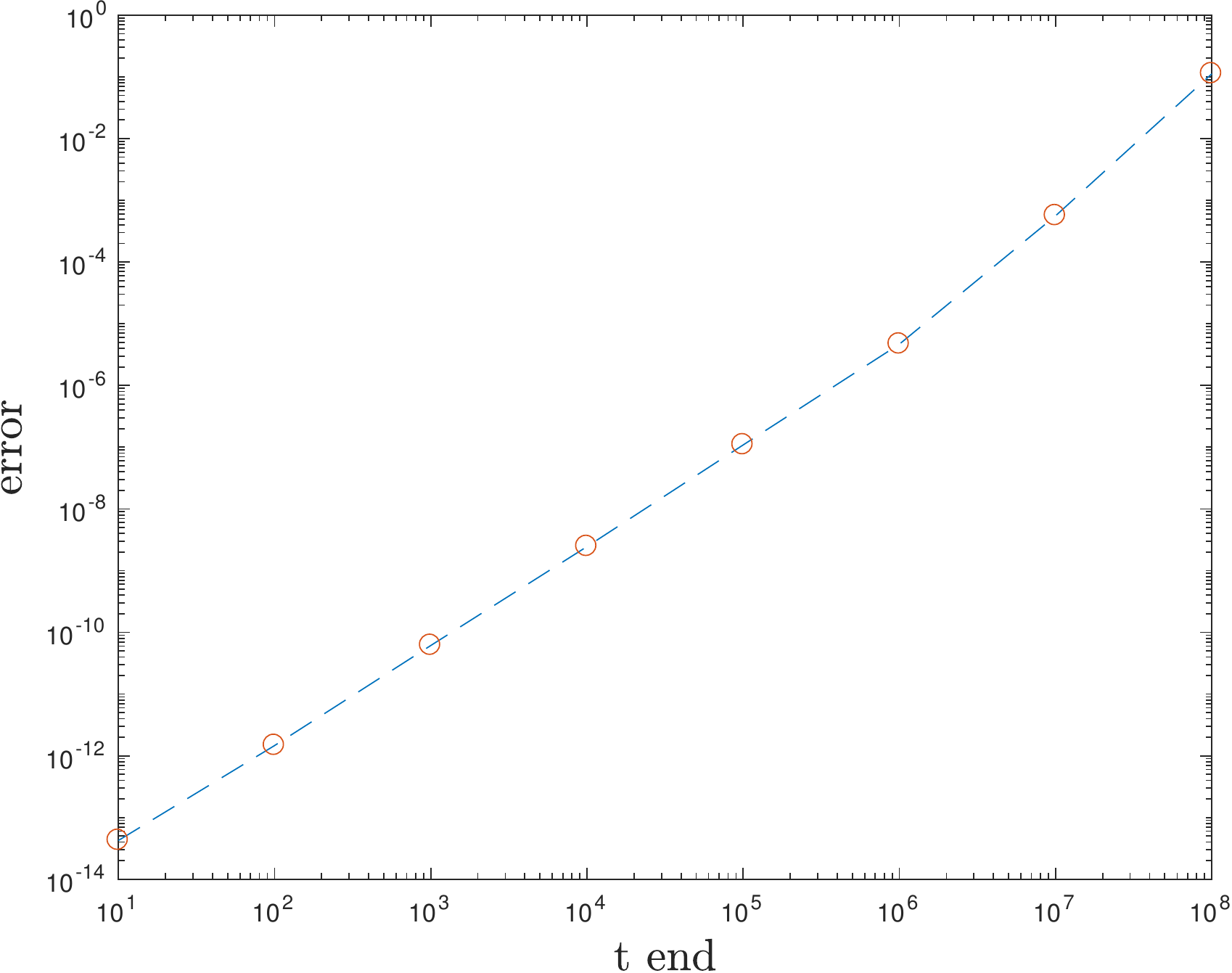}
  \caption{Log-log plot of divergence vs time. The planet model was integrated at two tolerances \protect\li{1e-13} and \protect\li{1e-15}. The norm of the difference between these solutions is plotted at $10^n$ TU, $n=1,\ldots,8$, showing no sign of chaotic behaviour for this set of ICs. \label{fg:planetdiverge}}
\end{figure}
\end{example}
\section{Conclusions and further work}

For two significant applications to do with DAEs, we have shown that differentiation of expressions, commonly done symbolically with the help of a computer algebra system, can be done efficiently and simply by AD.

\medskip
First, for the Dummy Derivatives index-reduction method a theoretical scheme is given that applies in principle to preparing a DAE for solution by any standard DAE or ODE initial value code.
DD-switching, which moves from one mode (local coordinate system) to another, is at the housekeeping level just a change of the size-\DOF set $S$ of indices $(j,l)$ that define the state vector.%

The scheme reduces order and index together, so one need not pre-reduce to first order form.
Finally, following this paper's theme, differentiating DAE components $f_i$ selectively (some more than others) does not need symbolic algebra; it can be done by standard AD methods of treating them as truncated power series.

We have proof-of-concept implementations in \matlab and \cpp.
It remains to be seen whether the scheme can be made efficient as a practical tool.
For DAEs from industrial applications that may need to switch among very many modes, it may (as in the more general case of hybrid systems) be worth keeping a run time data base of modes used, if this can speed up re-entry to a mode that has been met before.

\medskip
Second, we have shown that for a DAE, all or part of whose equations $f_i=0$ derive from the Lagrangian $L$ of a mechanical system, producing the $f_i$ from $L$ can be done by pure AD without symbolic algebra.
The theory was illustrated by simulation examples: a constrained mechanical system, a forced pendulum as a prescribed-path control problem, and an ODE of planetary motion.

The method of directly solving from a Lagrangian by overlaying one AD type on another might be used with other DAE solvers and AD tools.
However our infrastructure, of \daets with \fadbadpp and our Lagrangian facility has several advantages:
\begin{itemize}
  \item For any SA-friendly DAE, the user does not need to perform index/order reduction, since it is handled by \daets automatically.
	\item It avoids large symbolic expressions that a computer algebra system typically generates when converting to a form suitable for integration by a standard ODE/DAE solver.
	
  \item It can be programmed in a way that is intuitive and close to the mathematics.
  
  \item It gives remarkably fast code in the cases we have tried.
  \item Constrained ``first kind'' Lagrangians in cartesian coordinates are often simpler to formulate than unconstrained ``second kind'' ones in generalised coordinates.
For a high-index DAE solver such as \daets, possible obstacles posed by index reduction and DD-switching are absent, and constrained systems are as easy to solve as unconstrained, which may make ``first kind'' forms more attractive.
\end{itemize}

Future work will explore our Lagrangian approach on  a variety of research and engineering problems, and in particular rigid-body mechanics simulations and control problems. 
We are particularly interested in hybrid systems, because of their importance in industrial engineering.

\paragraph*{\bf Acknowledgments.} 
We acknowledge with thanks the support given to JDP by the Leverhulme Trust and the Engineering and Physical Sciences Research Council (EPSRC), both of the UK; 
and NN, GT, and XL by the Canadian Natural Sciences and Engineering Research Council (NSERC).
\newpage
\bibliographystyle{gOMS}

\appendix
\section{Extract from code for planetary problem}\label{sc:nbody3Drelative}

	This is the \li{fcn()} code for 3D motion of $n+1$ gravitating bodies, where body 0 is the Sun and the $x,y,z$ positions of the other bodies are relative to it.
	It was specialised to the problem in \exref{DETEST} by providing suitable input to the main program, not shown. Note this function is not restricted to 5 bodies: their number and masses are passed as parameters.

\medskip
\begin{lstlisting}
template <typename T>
void fcn( T t, const T *z, T *f, void *pp ) {
  const double *param = (double *)pp;
  const int  nMASS = param[0],
                 n = nMASS-1;
  const double   G = param[1]; 
  const double  *m = param + 2,
         *mplanet  = m+1;                // the masses EXCLUDING the Sun
  const double Mtotal = (m + nMASS)[0];  // total mass, calculated in main program

  typedef Vector3D< B<T> > vec3;
  vector< B<T> > q(3*n), qp(3*n);  // independent variables
  B<T> L;                          // for storing Lagrangian
  // { ... } ensures all intermediate variables go out of scope
  {
    init_q_qp(z,q,qp);  // setup q, qp

    // Convert to a vector of 3D vectors.
    vector< vec3 >  Q(n), Qp(n);
    for (int imass=0; imass<n; imass++) {
      Q[imass]  = vec3( q [3*imass],  q[3*imass+1], q [3*imass+2] );
      Qp[imass] = vec3( qp[3*imass], qp[3*imass+1], qp[3*imass+2] );
    }

    vector< vec3 > R(nMASS), Rp(nMASS);
    R[0] = Rp[0] = vec3(0,0,0);
    for (int imass=0; imass<n; imass++) {
      R [0] -= mplanet[imass]*Q [imass]; Rp[0] -= mplanet[imass]*Qp[imass];
    }
    R [0] /= Mtotal; Rp[0] /= Mtotal;
    
    // then set r_1, ..., r_n and their derivatives:
    for (int imass=1; imass<nMASS; imass++) {
      R [imass] = Q [imass-1]+R [0]; Rp[imass] = Qp[imass-1]+Rp[0];
    }
    
    // Compute KE and PE in terms of r and rp arrays
    B<T> KE = 0;
    for (int imass=0; imass<nMASS; imass++) 
      KE += m[imass] * Rp[imass]*Rp[imass];
    KE *= 0.5;
    
    // Potential Energy (sum of all mass-to-mass PEs, -> -oo as bodies
    // become close)
    B<T> PE = 0;
    for (int imass=0; imass<nMASS; imass++)
      for (int jmass=imass+1; jmass<nMASS; jmass++)
        PE -= m[imass] * m[jmass] / norm(R[imass]-R[jmass]);
    PE = G*PE; // bring in gravitational constant
    
    L = KE-PE;
  }
  vector<B<T> > C; // dummy constraint variable
  setupEquations(L,z,q,qp,C,f);
}
\end{lstlisting}
 
\end{document}